\documentclass[a4paper,11pt]{article}

\usepackage{amsmath}
\usepackage{amsthm}
\usepackage{amsfonts}

\newcommand{\mm}{\mathcal{M}}
\newcommand{\mq}{\mathbb{Q}}
\newcommand{\mr}{\mathbb{R}}
\newcommand{\mn}{\mathbb{N}}
\newcommand{\rP}{\mathbb{P}}
\newcommand{\rE}{\mathbb{E}}
\newcommand{\lix}{\underset{x\rightarrow\infty}{\lim}}

\newcommand{\su}{\underset{n\geq 0}\sup}
\newcommand{\sur}{\underset{n\geq 1}\sup}
\newcommand{\sun}{\underset{n\geq 0}\sup W_n}

\newtheorem{thm}{\noindent Theorem}[section]

\theoremstyle{remark}
\newtheorem{rem}{Remark}[section]
\newtheorem{ex}{Example}[section]

\begin{document}

\title{Regular variation in the branching random walk}\date{}
\author{Aleksander Iksanov\footnote {e-mail address:
iksan@unicyb.kiev.ua} and Sergey Polotskiy\footnote {e-mail
address: pilot$\_$ser@mail.ru} \\ \small{\emph{Faculty of
Cybernetics},
\emph{National T.Shevchenko University of Kiev}},\\
\small{\emph{01033 Kiev, Ukraine}}}\maketitle

\begin{abstract}
Let $\{\mm_n, n=0,1,\ldots\}$ be the supercritical branching
random walk starting with one initial ancestor located at the
origin of the real line. For $n=0,1,\ldots$ let $W_n$ be the
moment generating function of $\mm_n$ normalized by its mean.
Denote by $AW_n$ any of the following random variables: maximal
function, square function, $L_1$ and a.s. limit $W$, $\su
|W-W_n|$, $\su |W_{n+1}-W_n|$. Under mild moment restrictions and
the assumption that $\rP\{W_1>x\}$ regularly varies at $\infty$ it
is proved that $\rP\{AW_n>x\}$ regularly varies at $\infty$ with
the same exponent. All the proofs given are non-analytic in the
sense that these do not use Laplace-Stieltjes transforms. The
result on the tail behaviour of $W$ is established in two distinct
ways.


\end{abstract}

\section{An introduction, notation and results}

Let $\mm$ be a point process on $\mr$, i.e. random, locally finite
counting measure. Explicitly, $$\mm(A)(\omega):=
\sum_{i=1}^{J(\omega)} \delta_{X_i(\omega)}(A),$$ where $J:=\mm
(\mr)$, $\{X_i: i=\overline{1,J}\}$ are the points of $\mm$, $A$
is any Borel subset of $\mr$ and $\delta_x$ is the Dirac measure
concentrated at $x$. We assume that $\mm$ has no atom at $+
\infty$, and the $J$ may be deterministic or random, finite or
infinite with positive probability.

Let $\{\mm_n, n=0,1,\ldots\}$ be the \underline{b}ranching
\underline{r}andom \underline{w}alk (BRW), i.e. the sequence of
point processes which, for any Borel set $B\subseteq \mr$, are
defined as follows: $\mm_0(B)=\delta_0 (B)$,
\begin{equation*}
\mm_{n+1}(B):=\sum_r \mm_{n,r}(B-A_{n,r}), n=0,1,\ldots,
\end{equation*}
where $\{A_{n,r}\}$ are the points of $\mm_n$, and $\{\mm_{n,r}\}$
are independent copies of $\mm$. More detailed definition of the
BRW can be found in, for example, \cite{BiggKypr, Iks04, IksRos}.

In the case when $\rP\{J<\infty\}=1$ we assume that $\rE J>1$. In
the contrary case the condition holds automatically. Thus we only
consider the supercritical BRW. As a consequence,
$\rP\{\mm_n(\mr)>0 \ \  \text{for all} \  n\}>0$.

In what follows we use the notation that is generally accepted in
the literature on the BRW: $A_u$ denotes the position on $\mr$ of
a generic point $u=i_1\ldots i_n$; the record $|u|=n$ means that
the $u$ is a point of $\mm_n$; the symbol $\sum_{|u|=n}$ denotes
the summation over all points of $\mm_n$;
$\mathcal{F}_n=\sigma(\mm_1,\ldots, \mm_n)$ denotes the $\sigma
$-field generated by $\{\mm_k, k=1,\ldots, n\}$; $\mathcal{F}_0$
is the trivial $\sigma$-field.

Define the function
\begin{equation*}
m(y):=\rE \int_\mr e^{yx}\mm(dx)=\rE \sum_{|u|=1}e^{y A_u}\in
(0,\infty], y \in \mr,
\end{equation*}
and assume that there exists a $\gamma>0$ such that
$m(\gamma)<\infty$. Set $Y_u:=e^{\gamma A_u}/m^{|u|}(\gamma)$ and
\begin{equation*}
W_n:=m(\gamma)^{-n}\int_\mr e^{\gamma x}\mm_n(dx)=\sum_{|u|=n}
Y_u.
\end{equation*}
As is well-known (see, for example, \cite{Kingman}), the sequence
$\{(W_n, \mathcal{F}_n), n=0,1,\ldots\}$ is a non-negative
martingale. Notice that $W_0=\rE W_n=1$.

Let $\{d_n, n=1,2,\ldots\}$ be the martingale difference sequence,
i.e.
$$W_n=1+\sum_{k=1}^n d_k, n=1,2,\ldots$$ The square function $S$ and maximal function $W^\ast$ are
defined by $$S:=\left(1+\sum_{k=1}^\infty d_k^2\right)^{1/2} \ \
\text{and} \ \ W^\ast:=\sun. $$ Set also
$$S_n:=\left(1+\sum_{k=1}^n d_k^2\right)^{1/2}, n=1,2,\ldots \ \
\text{and} \ \ \Delta:=\sur |d_n|.$$ Recall that since $W_n$ is a
non-negative martingale all the defined variables are a.s. finite
(for finiteness of $S$ for general $L_1$-bounded martingales we
refer to \cite{Austin} or to Theorem 2 on p.390 \cite{Chow}).

When the martingale $W_n$ is uniformly integrable, we denote by
$W_\infty=W$ its $L_1$ and a.s. limit, and then define
$$M:=\su|W-W_n|=\su|\sum_{k=n+1}^\infty d_k|.$$

Lemma 1 \cite{IksNegad} (see also \cite{Bigg} for a slightly
different proof in the case $J<\infty$ a.s.) states that there
exist $r\in (0,1)$ and $\theta=\theta(r)>1$ such that whenever
$t>1$
\begin{equation}
\label{inn} \rP\{W>t\} \leq \rP \{W^\ast>t\} \leq \theta
\rP\{W>rt\}.
\end{equation}
This suggests that the tail behaviours of $W$ and $W^\ast$ are
quite similar.

Let now $\{f_n:=\sum_{k=0}^n g_k, n=0,1,\ldots\}$ be any
martingale. It is well-known that the distributions of maximal
$f^\ast:=\su |f_n|$ and square $S(f):=(\sum_{k=0}^\infty
g^2_k)^{1/2}$ functions are close in many respects. The evidence
in favor of such a statement is provided by, for example, the
(moment) Burkholder-Gundy-Davis inequality (Theorem 1.1
\cite{BuDa}) or the distribution function inequalities like
(\ref{bbb}) of this paper. From \cite{BuGu} and \cite{BuDa} and
many other subsequent works it follows that there exist a subset
$\mathcal{H}$ of the set of all martingales and a class
$\mathcal{A}$ of operators on martingales such that the
distributions of $A_1h$ and $A_2h$ are close in an appropriate
sense whenever $A_i\in \mathcal{A}$ and $h\in \mathcal{H}$. Often
can it be possible to express this closeness via moment or
distribution function inequalities like those mentioned above.
Keeping this in mind, it would not be an unrealistic conjecture
that the regular variation of $\rP\{A_1h>x\}$ is equivalent to
that of $\rP\{A_2h>x\}$, where $A_i$ and $h$ belong to some
subsets of operators and martingales respectively that may be
different from $\mathcal{A}$ and $\mathcal{H}$. On the other hand,
let us notice that as far as we know the conjecture does not
follow from previously known results on martingales.

The aim of this paper is to prove a variant of the conjecture for
the martingales $W_n$ and operators $A_i, i=\overline{1,5}$ given
as follows: $A_1W=W^\ast$, $A_2W=\Delta$, $A_3W=S$,
$A_4W=W_\infty$, $A_5W=M$.

In addition to the notation introduced above, other frequently
used notation and conventions include: $L(t)$ denotes a function
that slowly varies at infinity; $1_{A}$ denotes the indicator
function of the set $A$; $f(t) \sim g(t)$ is abbreviation of the
limit relation $\underset{t\rightarrow
\infty}{\lim}\dfrac{f(t)}{g(t)}=1$; $x^+:=\max(x,0)$; $x\wedge
y=\min (x,y)$; $x\vee y=\max (x,y)$; we write $\rP_n\{\cdot\}$
instead of $\rP\{\cdot|\mathcal{F}_n\}$, and $\rE_n\{\cdot\}$
instead of $\rE\{\cdot|\mathcal{F}_n\}$; the record "const"
denotes a constant whose value is of no importance and may be
different on different appearances.

Now we are ready to state our result.
\begin{thm}
\label{t1} Assume that there exist $\beta>1$ and $\epsilon>0$ such
that
\begin{equation}
\label{mom} k_\beta:=\rE \sum_{|u|=1}Y_u^\beta<1,\ \ \rE
\sum_{|u|=1}Y_u^{\beta+\epsilon}<\infty \ \ \text{and}
\end{equation}
\begin{equation}
\label{var}\rP\{W_1>x\}\sim x^{-\beta}L(x).
\end{equation}
Then
\newline (I) $\rP\{W^\ast>x\}\sim
\rP\{\Delta>x\}\sim \rP\{S>x\}\sim
(1-k_\beta)^{-1}\rP\{W_1>x\}$;\newline (II) $W_n$ converges almost
surely and in mean to a random variable $W$ and
\begin{equation}
\label{ww} \rP\{W>x\}\sim (1-k_\beta)^{-1}\rP\{W_1>x\};
\end{equation}
$$\rP\{M>x\}\sim (1-k_\beta)^{-1}\rP\{W_1>x\}.$$
\end{thm}
\begin{rem}
We are not aware of any papers on branching processes which
investigate the tail behaviour of random variables like $\Delta$,
$M$ or $S$. \cite{IksNegad} is the only paper we know of that
deals with the tail behaviour of random variables like $W^\ast$.
\end{rem}
\begin{rem}
When $\gamma=0$ and $J<\infty$ a.s., $W_n$ reduces to the
(supercritical) normalized Galton-Watson process. In this case
(\ref{ww}) was proved in \cite{BiDo74} for non-integer $\beta$ and
in \cite{Meyer82} for integer $\beta$. When $\gamma>0$, $J<\infty$
a.s. and $\mm(-\infty, -\gamma^{-1}\log m(\gamma))=0$ a.s., $W$
can be viewed as a limit random variable in the Crump-Mode
branching process. In this case (\ref{ww}) was established in
\cite{BiDo75} for non-integer $\beta$. The technique used in the
last three cited works is purely analytic (based on using the
Laplace-Stieltjes transforms and Abel-Tauberian theorems) and
completely different from ours. On the other hand, let us notice
that the above mentioned analytic approach was successfully
employed and further developed by the second-named author. In
2003, in an unpublished diploma paper he proved (\ref{ww}) for
non-integer $\beta$ for the general case treated here.

Our desire to find a non-analytic proof of (\ref{ww}) was a
starting point for the development of this paper. In the course of
writing two different (non-analytic) proofs were found. One of
these proofs given in Section 2 falls within the general scope of
the paper. The second given in Section 3 continues a line of
research initiated in \cite{Iks04}, \cite{IksRos}, \cite{Iks06}.
Here an underlying idea is that the martingale $W_n$ and so called
\emph{perpetuities} have many features in common. In particular,
several non-trivial results on perpetuities (however, it seems,
only those related to perpetuities with not all moments finite)
can be effectively exploited to obtain similar results on the
limiting behaviour of $W_n$. Maybe we should recall that, in
modern probability, by a perpetuity is meant a random variable
$$B_1+\sum_{k=2}^\infty A_1A_2\cdots A_{k-1}B_k,$$ provided the latter series absolutely converges, and where
$\{(A_k, B_k): k=1,2,\ldots\}$ are independent identically
distributed random vectors.
\end{rem}

The paper is structured as follows. In Section 2 we prove Theorem
\ref{t1}. Here an essential observation is that, given
$\mathcal{F}_n$, $W_{n+1}$ looks like a weighted sum of
independent identically distributed random variables. This allows
us to exploit the well-known result \cite{MikSam} on the tail
behaviour of such sums under the regular variation assumption. The
second key ingredient of the proof is using the distribution
function inequalities for martingales. In Section 3 we give
another proof of (\ref{ww}) which rests on a relation between the
BRW and perpetuities. Here availability of Grincevi\v{c}ius-Grey
\cite{Grey} result on the tail behaviour of perpetuities is
crucial. Finally, in Section 4 we discuss applicability of Theorem
\ref{t1} to several classes of point processes. The section closes
with two remarks which show that (\ref{mom}) and (\ref{var}) are
not necessary conditions for the regular variation of the tails of
$W^\ast$, $W$ and a related random variable.

\section{Proof of Theorem \ref{t1}}

(I) We will prove the result for $W^\ast$ and $\Delta$
simultaneously. To this end, let $Q$ and $\tilde{Q}$ be
independent identically distributed random variables whose
distribution is supported by $(a,\infty), a>-\infty$. Assume that
$\rP\{Q>x\}\sim x^{-\beta}L(x)$ for $\beta>1$. In particular, this
assumption ensures that $\rE |Q|<\infty$ and $\rP\{|Q|>x\}\sim
\rP\{Q>x\}$. With a slight abuse of notation, set
$Q^s:=|Q|-|\tilde{Q}|$. Then
\begin{equation}
\label{sym} 1-F(x):=\rP\{|Q^s|>x\}\sim 2x^{-\beta}L(x).
\end{equation}
Indeed, $1-F(x)=2\int_0^\infty (1-G(x+y))dG(y)$, where
$G(x)=\rP\{|Q|\leq x\}, x\geq 0$. Now (\ref{sym}) follows from
monotonicity of $1-G$, the relation $1-G(x+y)\sim 1-G(x), y\in
\mr$ and Fatou's lemma.

The equality
\begin{equation*}
\rE t(Z)=\rE \sum_{|u|=1}Y_ut(Y_u),
\end{equation*}
which is assumed to hold for all bounded Borel functions $t$,
defines the distribution of a random variable $Z$. More generally,
\begin{equation}
\label{ex} \rE t(Z_1\cdots Z_n)=\rE \sum_{|u|=n}Y_u t(Y_u),
\end{equation}
where $Z_1, Z_2, \ldots$ are independent copies of the $Z$. Notice
that we can permit for (\ref{ex}) to hold for any Borel function
$t$. In that case we assume that if the right-hand side is
infinite or does not exist, the same is true for the left-hand
side.

Under the assumptions of the theorem, the function $k_x:=\rE
\sum_{|u|=1}Y_u^x$ is log-convex for $x\in (1,\beta)$, $k_1=1$ and
$k_\beta<1$. Therefore,
\begin{equation}
\label{innn} k_{\beta-\epsilon}<1 \ \ \text{for all} \ \epsilon
\in (0,\beta-1).
\end{equation}
Also we can pick a $\delta\in (0,\beta-1)$ such that
$k_{\beta+\delta}<1$. By using these facts and equality (\ref{ex})
we conclude that with this $\delta$
\begin{equation}
\label{useful} \rE
\sum_{|u|=n}Y_u^{\beta-\delta}=k^n_{\beta-\delta}<1 \ \ \text{and}
\ \ \rE\sum_{|u|=n}Y_u^{\beta+\delta}=k^n_{\beta+\delta}<1.
\end{equation}
Let us notice, for later needs, that we can choose $\delta$ as
small as needed. Among other things, (\ref{useful}) implies that
for $x\in [1,\beta+\delta]$
\begin{equation}
\label{inter} \sum_{|u|=n}Y_u^x<\infty \ \ \text{a.s.}
\end{equation}

Until further notice, we fix an arbitrary $n\in \mn$. Put
$$T_n:=|\sum_{|u|=n}Y_u Q_u| \ \ \text{and} \ \ X_n:=\sum_{|u|=n}Y_u |Q_u|.$$
Given $\mathcal{F}_n$, let $\{Q_u:|u|=n\}$ and $\{Q_u^s:|u|=n\}$
be conditionally independent copies of the random variables $Q$
and $Q^s$ respectively. In view of (\ref{inter}), an appeal to
Lemma A3.7\cite{MikSam} allows us to conclude that
\begin{equation}
\label{basic} \rP_n\{T_n>x\}\sim \sum_{|u|=n}Y_u^\beta
\rP\{|Q|>x\} \ \ \text{a.s.}
\end{equation}
The cited lemma assumes that each term of the series on the
left-hand side has zero mean, but this condition is not needed in
the proof of the result used here.

Denote by $\mu_n^{\mathcal{F}_n}$ the conditional on
$\mathcal{F}_n$ median of $X_n$, i.e. $\mu_n^{\mathcal{F}_n}$ is a
random variable that satisfies
$$\rP_n\{X_n-\mu_n^{\mathcal{F}_n}\geq 0\}\geq 1/2 \leq \rP_n\{X_n-\mu_n^{\mathcal{F}_n}\leq 0\} \ \
\text{a.s.}$$ Let also $\mu_n$ denote the usual median of $X_n$.
Since $\mu_n^{\mathcal{F}_n}\geq 0$ a.s., we have from
(\ref{basic})
$$\underset{x\to\infty}{\lim \sup} \dfrac{\rP_n\{T_n>x+\mu_n^{\mathcal{F}_n}\}}{\rP\{|Q|>x\}}\leq \sum_{|u|=n}Y_u^\beta \ \
\text{a.s.}$$ If we could prove that for large $x$
\begin{equation}
\label{ineq}
\dfrac{\rP_n\{T_n>x+\mu_n^{\mathcal{F}_n}\}}{{\rP\{|Q|>x\}}}\leq
U_n \ \ \text{a.s. \ \ and} \ \ \rE U_n<\infty,
\end{equation}
where $U_n$ is a random variable, then using Fatou's lemma yielded
\begin{equation}
\label{inter2} \underset{x\to\infty}{\lim \sup}\rE
\dfrac{\rP_n\{T_n>x+\mu_n^{\mathcal{F}_n}\}}{\rP\{|Q|>x\}}=
\underset{x\to\infty}{\lim \sup}
\dfrac{\rP\{T_n>x+\mu_n\}}{\rP\{|Q|>x\}}\leq \rE
\sum_{|u|=n}Y_u^\beta\overset{(\ref{ex})}{=}k^n_\beta.
\end{equation}
Since $\rP\{|Q|>x+\mu_n\}\sim \rP\{|Q|>x\}$, (\ref{inter2})
implied that
$$\underset{x\to\infty}{\lim \sup}
\dfrac{\rP\{T_n>x\}}{\rP\{|Q|>x\}}\leq k^n_\beta.$$ On the other
hand, by using (\ref{basic}) and Fatou's lemma the reverse
inequality for the lower limit follows easily. Therefore, as soon
as (\ref{ineq}) is established, we get
\begin{equation}
\label{main} \rP\{T_n>x\}\sim k^n_\beta \rP\{|Q|>x\}.
\end{equation}

We now intend to show that (\ref{ineq}) holds with
\begin{equation}
\label{U} U_n=const
\left(\sum_{|u|=n}Y_u^{\beta-\delta}+\sum_{|u|=n}Y_u^{\beta+\delta}
\right),
\end{equation}
for appropriate small $\delta$ that satisfies (\ref{useful}).
Notice that
\begin{equation}
\label{expec} \rE U_n=const
(k^n_{\beta-\delta}+k^n_{\beta+\delta})<\infty.
\end{equation}

By the triangle inequality and conditional symmetrization
inequality
\begin{equation}
\label{start} (1/2)\rP_n\{T_n>x+\mu_n^{\mathcal{F}_n}\}\leq
(1/2)\rP_n\{X_n>x+\mu_n^{\mathcal{F}_n}\}\leq
\rP_n\{|\sum_{|u|=n}Y_uQ_u^s|>x\}.
\end{equation}
Let us show that for $x>0$
$$\rP_n\{|\sum_{|u|=n}Y_uQ_u^s|>x\}\leq$$
\begin{equation}
\label{imp}\leq
\rP_n\{\underset{|u|=n}{\sup}Y_u|Q_u^s|>x\}+x^{-2}\rE_n\left(\sum_{|u|=n}Y_u^2(Q_u^s)^21_{\{Y_u|Q_u^s|\leq
x\}}\right):=I_1(n,x)+I_2(n,x).
\end{equation}
Let $\{Y(k)Q^s(k):k=1,2,\ldots\}$ be any enumeration of the set
$\{Y_uQ_u^s:|u|=n\}$. The inequality $\rE |Q|<\infty$ implies that
the series $\sum_{|u|=n}Y_uQ_u$ is absolutely convergent.
Therefore $\sum_{|u|=n}Y_uQ_u=\sum_{k=1}^\infty Y(k)Q^s(k)$.
Define
\begin{equation*}
\tau_x:=
\begin{cases}
\inf\{k\geq 1: Y(k)|Q^s(k)|>x\}, \text{\ \ if \ \ }
\underset{k\geq
1}{\sup}Y(k)|Q^s(k)|>x; \\
+\infty, \text{\ \ otherwise }.
\end{cases}
\end{equation*}
For any fixed $m\in \mn$ and $x>0$
$$\rP_n\{|\sum_{k=1}^m Y(k)Q^s(k)|>x\}\leq$$
$$\leq \rP_n\{\tau_x\leq m-1\}+\rP_n\{|\sum_{k=1}^m Y(k)Q^s(k)|>x,
\tau_x\geq m\}\leq$$$$\leq \rP_n\{\underset{1\leq k\leq
m-1}{\sup}Y(k)|Q^s(k)|>x\}+\rP_n\{|\sum_{k=1}^{\tau_x\wedge
m}Y(k)Q^s(k)|>x\}\leq$$ (by Markov inequality) $$\leq
\rP_n\{\underset{1\leq k\leq
m-1}{\sup}Y(k)|Q^s(k)|>x\}+x^{-2}\rE_n\left(\sum_{k=1}^m
Y(k)Q^s(k)1_{\{\tau_x\geq k\}}\right)^2\leq $$ ($\rE_n Q^s(k)=0$
and, given $\mathcal{F}_n$, $Q^s(k)$ and $1_{\{\tau_x\geq k\}}$
are independent) $$\leq \rP_n\{\underset{1\leq k\leq
m-1}{\sup}Y(k)|Q^s(k)|>x\}+x^{-2}\rE_n \sum_{k=1}^m
Y^2(k)(Q^s(k))^21_{\{Y(k)|Q^s(k)|\leq x\}}.$$ If the distribution
of $Q^s$ is continuous, sending $m\to\infty$ then completes the
proof of (\ref{imp}).

Assume now that the distribution of $Q^s$ has atoms. Let $R$ be a
random variable with uniform distribution on $[-1,1]$, which is
independent of $Q^s$. Given $\mathcal{F}_n$, let $\{R_u:|u|=n\}$
be conditionally independent copies of $R$ which are also
independent of $\{Q_u^s:|u|=n\}$. Since for all $t>0$
$$\rP\{|Q^s|>t\}\leq 2\rP\{|Q^s||R|>t/2\},$$ we have by Theorem
3.2.1\cite{Woy}
\begin{equation}
\label{cont} \rP_n\{|\sum_{|u|=n}Y_uQ_u^s|>t\}\leq
4\rP_n\{|\sum_{|u|=n}Y_uQ_u^sR_u|>t/4\},
\end{equation}
and the distribution of $Q^sR$ is (absolutely) continuous. Now we
can apply the already established part of (\ref{imp}) to the
right-hand side of (\ref{cont}). Strictly speaking, when the
distribution of $Q^s$ has atoms, (\ref{imp}) should be written in
a modified form: additional constants should be added, and $Q_u^s$
should be replaced with $Q_u^sR_u$. On the other hand, a perusal
of the subsequent proof reveals that only the regular variation of
$\rP\{|Q^s|>x\}$ plays a crucial role. Therefore, for ease of
notation we prefer to keep (\ref{imp}) in its present form. This
does not cause any mistakes as $\rP\{|Q^sR|>x\}\sim \rE |R|^\beta
\rP\{|Q^s|>x\}$.

Assume temporarily that $1-F(x)$ regularly varies with index
$-\beta$, $\beta\in (1,2)$. Set $T(x):=\int_0^x y^2dF(y)$. By
Theorem 1.6.4 \cite{BGT}
\begin{equation*}
T(x)\sim \dfrac{\beta}{2-\beta} x^2(1-F(x))\sim
\dfrac{\beta}{2-\beta}x^{2-\beta}L_1(x).
\end{equation*}
Also by Theorem 1.5.3 \cite{BGT} there exists a non-decreasing
$S(x)$ such that
\begin{equation}
\label{s} T(x)\sim S(x).
\end{equation}
For any $A_i>0$ and $\delta$ defined in (\ref{useful}) there
exists an $x_i>0$ such that whenever $x\geq x_i, i=1,2,3$
\begin{equation}
\label{A1} x^{\beta+\delta}(1-F(x))\geq 1/A_1;
\end{equation}
\begin{equation}
\label{A2} x^{\beta-2+\delta}S(x)\geq 1/A_2;
\end{equation}
\begin{equation}
\label{A3} T(x)\leq
(A_3+\dfrac{\beta}{2-\beta})x^2(1-F(x)):=Bx^2(1-F(x)).
\end{equation}
Also for any $A_i>1$ and the same $\delta$ as above there exists
an $x_i>0$ such that whenever $x\geq x_i$ and $ux\geq x_i,
i=4,5,6$
\begin{equation}
\label{A4} \dfrac{1-F(ux)}{1-F(x)}\leq A_4 (u^{-\beta+\delta}\vee
u^{-\beta-\delta});
\end{equation}
\begin{equation}
\label{A5} \dfrac{T(ux)}{T(x)}\leq A_5 (u^{2-\beta+\delta}\vee
u^{2-\beta-\delta});
\end{equation}
\begin{equation}
\label{A6} \dfrac{T(ux)}{T(x)}\leq A_6\dfrac{S(ux)}{S(x)}.
\end{equation}
(\ref{A4}) and (\ref{A5}) follows from Potter's bound (Theorem
1.5.6 (iii)\cite{BGT}); (\ref{A6}) is implied by (\ref{s}). Set
$x_0:=\underset{1\leq i \leq 6}{\max}x_i$ and assume that $x_0>1$.

To check (\ref{ineq}) and (\ref{U}), we consider three cases: (a)
$\beta\in (1,2)$; (b) $\beta>2, \beta\neq 2^n, n\in\mn$; (c)
$\beta=2^n, n\in \mn$.

(a)  For any fixed $x\geq x_0$ and $y>0$
$$\dfrac{I_1(n,x/y)}{1-F(x)} \leq \sum_{|u|=n}\dfrac{\rP_n\{Y_u|Q_u^s|>x/y\}}{1-F(x)}=
\sum_{|u|=n}\dfrac{1-F(x(yY_u)^{-1})}{1-F(x)}=$$
$$=\sum_{|u|=n}\cdots 1_{\{yY_u>x/x_0\}}+\sum_{|u|=n}\cdots 1_{\{yY_u\leq
x/x_0\}}=:I_{11}(n,x,y)+I_{12}(n,x,y).$$ Since
$$(yY_u)^{\beta+\delta}\geq
(yY_u)^{\beta+\delta}1_{\{yY_u>x/x_0\}}\geq
(x/x_0)^{\beta+\delta}1_{\{yY_u>x/x_0\}},$$ we get
$$I_{11}(n,x,y)\leq
(x_0y)^{\beta+\delta}\dfrac{\sum_{|u|=n}Y_u^{\beta+\delta}}{x^{\beta+\delta}(1-F(x))}\overset{(\ref{A1})}{\leq}
A_1x_0^{\beta+\delta}y^{\beta+\delta}\sum_{|u|=n}Y_u^{\beta+\delta}.$$
Further
$$I_{12}(n,x,y)\overset{(\ref{A4})}{\leq}A_4\sum_{|u|=n}(yY_u)^{\beta-\delta}\vee(yY_u)^{\beta+\delta}\leq
A_4\left(y^{\beta-\delta}\sum_{|u|=n}Y_u^{\beta-\delta}+y^{\beta+\delta}\sum_{|u|=n}Y_u^{\beta+\delta}\right);
$$ $$\dfrac{I_2(n,x/y)}{1-F(x)}=y^2\sum_{|u|=n}\dfrac{Y^2_u\int_0^{x(yY_u)^{-1}}z^2dF(z)}{x^2(1-F(x))}\overset{(\ref{A3})}{\leq}
By^2\sum_{|u|=n}\dfrac{Y^2_uT(x(yY_u)^{-1})}{T(x)}=
$$$$=By^2\left(\sum_{|u|=n}\cdots 1_{\{yY_u>
x/x_0\}}+\sum_{|u|=n}\cdots 1_{\{yY_u\leq
x/x_0\}}\right)=:By^2(I_{21}(n,x,y)+I_{22}(n,x,y)).$$
$$I_{21}(n,x,y)\overset{(\ref{A6})}{\leq}A_6\sum_{|u|=n}\dfrac{Y^2_uS(x(yY_u)^{-1})}{S(x)}1_{\{yY_u>
x/x_0\}}\leq $$ $$\leq
A_6(x_0y)^{\beta-2+\delta}S(x_0)\sum_{|u|=n}\dfrac{Y_u^{\beta+\delta}}{x^{\beta-2+\delta}S(x)}\overset{(\ref{A2})}{\leq}
A_2A_6(x_0y)^{\beta-2+\delta}S(x_0)\sum_{|u|=n}Y_u^{\beta+\delta};$$
$$I_{22}(n,x,y)\overset{(\ref{A5})}{\leq} A_5 \sum_{|u|=n}Y_u^2((yY_u)^{\beta-2-\delta}\vee
(yY_u)^{\beta-2+\delta})\leq$$
$$\leq A_5\left(y^{\beta-2-\delta}\sum_{|u|=n}Y_u^{\beta-\delta}+y^{\beta-2+\delta}\sum_{|u|=n}Y_u^{\beta+\delta}\right).$$
Thus according to (\ref{imp}) we have proved that for $x\geq x_0$
and $y>0$
\begin{equation}
\label{impo}
\dfrac{\rP_n\{|\sum_{|u|=n}Y_uQ_u^s|>x/y\}}{\rP\{|Q^s|>x\}}\leq
const
\left(y^{\beta-\delta}\sum_{|u|=n}Y_u^{\beta-\delta}+y^{\beta+\delta}\sum_{|u|=n}Y_u^{\beta+\delta}
\right).
\end{equation}
In particular, since $\rP\{|Q^s|>x\}\sim 2\rP\{|Q|>x\}$ then
setting in (\ref{impo}) $y=1$ and using (\ref{start}) leads to
(\ref{ineq}) with $U_n$ being a multiple of the right-hand side of
(\ref{impo}).

In the remaining cases we only investigate the situation when
$\beta\in (2,4)$ and $\beta=2$. For other values of $\beta$ the
inequality (\ref{ineq}) with $U_n$ satisfying (\ref{U}) follows by
induction.

(b) $\beta\in (2,4)$. Given $\mathcal{F}_n$, let $\{\tilde{N},
N_u:|u|=n\}$ be conditionally independent copies of a random
variable $N$ with normal $(0,1)$ distribution. Using the approach
similar to that exploited to obtain (\ref{cont})(this fruitful
argument has come to our attention from \cite{MikSam}) allows us
to conclude that for $x>0$ and appropriate positive constants
$c_1, c_2$
$$\rP_n\{|\sum_{|u|=n}Y_uQ_u^s|>x\}\leq$$
\begin{equation}
\label{sam} \leq
c_1\rP_n\{|\sum_{|u|=n}Y_uN_uQ_u^s|>c_2x\}=c_1\rP_n\left\{|\tilde{N}|\left(\sum_{|u|=n}Y^2_u(Q_u^s)^2\right)^{1/2}>c_2x\right\}.
\end{equation}
Notice that $\rP\{(Q^s)^2>x\}$ regularly varies with index
$-\beta/2\in (-2,-1)$. Also it is obvious that, if needed, we can
reduce $\delta$ in (\ref{useful}) to ensure that
$k_{\beta-2\delta}<1$ and $k_{\beta+2\delta}<1$. Therefore we can
use (\ref{impo}) with $Y_u$ replaced with $Y_u^2$, and $Q_u^s$
replaced with $(Q_u^s)^2$ which gives after a little manipulation
that for $x\geq x_0^{1/2}$
$$\dfrac{\rP_n\{|\sum_{|u|=n}Y_uQ_u^s|>x\}}{\rP\{|Q^s|>x\}}\leq const\left(\sum_{|u|=n}Y_u^{\beta+2\delta}\rE|N|^{\beta+2\delta}+
\sum_{|u|=n}Y_u^{\beta-2\delta}\rE|N|^{\beta-2\delta}\right).$$ By
using the same argument as in the case (a) we can check that
(\ref{ineq}) and (\ref{U}) have been proved.

(c) $\beta=2$. In the same manner as we have established
(\ref{impo}) it can be proved that for $y>0$ and large $x$
$$\dfrac{\rP_n\{\sum_{|u|=n}Y^2_u(Q_u^s)^2>x/y\}}{\rP\{(Q^s)^2>x\}}\leq
const
\left(y^{2-2\delta}\sum_{|u|=n}Y_u^{2-2\delta}+y^{2+2\delta}\sum_{|u|=n}Y_u^{2+2\delta}
\right).$$ Hence an appeal to (\ref{sam}) assures that
(\ref{ineq}) and (\ref{U}) hold in this case too.

We have a representation
\begin{equation}
\label{repr} W_{n+1}=\sum_{|u|=n}Y_uW_1^{(u)},
\end{equation}
where, given $\mathcal{F}_n$, $W_1^{(u)}$ are (conditionally)
independent copies of $W_1$. Each element of the set
$\{W_1^{(u)}:|u|=n\}$ is constructed in the same way as $W_1$, the
only exception being that while $W_1$ is defined on the whole
family tree, $W_1^{(u)}$ is defined on the subtree with root $u$.

We only give a complete proof for the $\Delta$. The analysis of
the $W^\ast$ is similar but simpler, and hence omitted. From
(\ref{repr}) we conclude that $|d_{n+1}|$ is the same as $T_n$
with $Q_u=W_1^{(u)}-1$. Hence by (\ref{basic})
\begin{equation}
\label{co}1_{\{\underset{1\leq k\leq n}{\max}|d_k| \leq
x\}}\rP_n\{|d_{n+1}|>x\} \sim \sum_{|u|=n}Y_u^\beta
\rP\{|W_1-1|>x\}\sim \sum_{|u|=n}Y_u^\beta \rP\{W_1>x\}  \ \
\text{a.s.},
\end{equation}
and by (\ref{main})
\begin{equation}
\label{nnn} \rP\{|d_{n+1}|>x\}\sim k^n_\beta \rP\{|W_1-1|>x\}\sim
k^n_\beta \rP\{W_1>x\}.
\end{equation}
Recall that
$$\rP\{\Delta>x\}=\rP\{|d_1|>x\}+\sum_{n=1}^\infty \rP\{\underset{1\leq k\leq n}{\max}|d_k| \leq
x,|d_{n+1}|>x\}=$$ $$=\rP\{|d_1|>x\}+\rE \sum_{n=1}^\infty
1_{\{\underset{1\leq k\leq n}{\max}|d_k| \leq
x\}}\rP_n\{|d_{n+1}|>x\}.$$

Using this, (\ref{co}) and applying Fatou's lemma twice allows us
to conclude that
\begin{equation*}
\underset{x\to\infty}{\lim \inf} \dfrac{\rP\{\Delta>x\}}{\rP
\{W_1>x\}}\geq 1+\sum_{n=1}^\infty \rE \underset{x\to\infty}{\lim
\inf} \dfrac{1_{\{\underset{1\leq k\leq n}{\max}|d_k| \leq
x\}}\rP_n\{|d_{n+1}|>x\}}{\rP\{W_1>x\}}\geq
\end{equation*}
\begin{equation*}
\geq 1+\sum_{n=1}^\infty \rE\left( \sum_{|u|=n}
Y_u^\beta\right)=(1-k_{\beta})^{-1}.
\end{equation*}

To complete the proof for $\Delta$ we must calculate the
corresponding upper limit. For this it suffices to check that for
large $x$ and large $n\in \mn$
\begin{equation}
\label{qqq} \dfrac{\rP\{|d_{n+1}|>x\}}{\rP\{W_1>x\}}\leq C_n \ \
\text{and} \ \ C_n \ \ \text{is a summable sequence},
\end{equation}
and use the dominated convergence theorem. Taking the expectation
in (\ref{ineq}) allows us to conclude that for $n=1,2,\ldots$ and
large $x$
$$\dfrac{\rP\{|d_{n+1}|>x+\mu_n\}}{{\rP\{W_1>x\}}}\leq const \rE U_n,$$
where $\mu_n$ is the median of
$V_n:=\sum_{|u|=n}Y_u|W_1^{(u)}-1|$, and $\rE U_n$ is given by
(\ref{expec}). The family of distributions of $V_n$ is tight. In
view of (\ref{nnn}), $$\rP\{|d_{n+1}|>x+y\}\sim \rP\{|d_{n+1}|>x\}
\ \ \text{locally uniformly in} \ y .$$ Therefore, (\ref{qqq})
holds with $C_n=const \rE U_n$ and the result for $\Delta$ has
been proved.

For later needs let us notice here that in the same way as above
we can prove that for fixed $n\in \mn$
\begin{equation}
\label{qq} \rP\{\underset{m\geq n}{\sup}W_m> x\}\sim
k_\beta^n(1-k_\beta)^{-1}\rP\{W_1>x\}.
\end{equation}

Consider now the square function $S$. Since $S\geq \Delta$ a.s.,
and we have already proved that $\rP\{\Delta>x\} \sim
(1-k_\beta)^{-1}\rP\{W_1>x\}$ then
$$\underset{x\to\infty}{\lim
\inf}\dfrac{\rP\{S>x\}}{\rP\{W_1>x\}}\geq \dfrac{1}{1-k_\beta}.$$
Therefore we must only calculate the upper limit. We begin with
showing that for any $n\in \mn$
\begin{equation}
\label{ss} \underset{x\to\infty}{\lim
\sup}\dfrac{\rP\{S_n>x\}}{\rP\{W_1>x\}}\leq
\sum_{m=0}^{n-1}k_\beta^m.
\end{equation}
We use induction on $n$.\newline (1) If $n=1$ then $S_1\leq W_1$
and (\ref{ss}) is obvious.\newline (2) Assume that (\ref{ss})
holds for $n=j$ and show that it holds for $n=j+1$. For every
$x>0$ and $\epsilon \in (0,1)$
$$\rP\{S_{j+1}>x\}\leq \rP\{S_j^2>(1-\epsilon)x^2\}+\rP\{d_{j+1}^2>(1-\epsilon)x^2\}+\rP\{S_j^2>\epsilon x^2, d_{j+1}^2
>\epsilon x^2\}=$$ write the latter $\rP$ as $\rE\rP_j$ and use $\mathcal{F}_j$-measurability of $S_j$ to get
$$=\rP\{S_j>(1-\epsilon)^{1/2}x\}+\rP\{|d_{j+1}|>(1-\epsilon)^{1/2}x\}+\rE
1_{\{S_j>\epsilon^{1/2}x\}}\rP_j\{|d_{j+1}|>\epsilon^{1/2}x\}.$$
According to (\ref{basic}) with $Q_u$ replaced by $W_1^{(u)}-1$,
$$\lix 1_{\{S_j>\epsilon^{1/2}x\}}\dfrac{\rP_j\{|d_{j+1}|>\epsilon^{1/2}x\}}{\rP\{W_1>x\}}=0 \ \
\text{a.s.},$$ and there exists a $\delta_1>0$ such that for large
$x$
$$1_{\{S_j>\epsilon^{1/2}x\}}\dfrac{\rP_j\{|d_{j+1}|>\epsilon^{1/2}x\}}{\rP\{W_1>x\}}\leq
\epsilon^{-\beta/2}\sum_{|u|=n}Y_u^\beta+\delta_1 \ \
\text{a.s.}$$ Therefore, by the dominated convergence $$\lix \rE
1_{\{S_j>\epsilon^{1/2}x\}}\dfrac{\rP_j\{|d_{j+1}|>\epsilon^{1/2}x\}}{\rP\{W_1>x\}}=0.$$
By the inductive assumption and (\ref{nnn})
$$\underset{x\to\infty}{\lim
\sup}\dfrac{\rP\{S_{j+1}>x\}}{\rP\{W_1>x\}}\leq
(1-\epsilon)^{-\beta/2}\sum_{m=0}^{j-1}k_\beta^m+(1-\epsilon)^{-\beta/2}k_\beta^j=(1-\epsilon)^{-\beta/2}\sum_{m=0}^jk_\beta^m.$$
Sending $\epsilon\to 0$ proves (\ref{ss}).

For $m=0,1,\ldots$ and fixed $n\in \mn$ set $\tilde{W}_m:=W_{m\vee
n}$ and $\tilde{\mathcal{F}}_m;=\mathcal{F}_{m\vee n}$. Choose
$\rho\in (0,\sqrt{3})$ so small that
$\nu:=\dfrac{2\rho^2}{3-\rho^2}2^{\beta+1}\in (0,1)$. Applying
Theorem 18.2 \cite{Burk} (in the notation of that paper take
$\beta=2$ and $\delta=\rho$) to the non-negative martingale
$(\tilde{W}_m, \tilde{\mathcal{F}}_m)$ gives
$$\rP\{(\sum_{m=n+1}^\infty d_m^2)^{1/2}>2x\}\leq \rP\{\underset{m\geq n}{\sup}\tilde{W}_m>\rho x\}+
\rP\{(\sum_{m=n+1}^\infty d_m^2)^{1/2}>2x, \ \underset{m\geq
n}{\sup}\tilde{W}_m\leq \rho x\}\leq $$
\begin{equation}
\label{bbb} \leq \rP\{\underset{m\geq n}{\sup}W_m>\rho
x\}+\dfrac{2\rho^2}{3-\rho^2}\rP\{(\sum_{m=n+1}^\infty
d_m^2)^{1/2}>x\}.
\end{equation}
By Potter's bound we can take $y>0$ such that
$\dfrac{\rP\{W_1>x\}}{\rP\{W_1>2x\}}\leq 2^{\beta+1}$ for $x\geq
y$. Set $A(y):=\underset{x\geq y}{\sup}\dfrac{\rP\{\underset{m\geq
n}{\sup}W_m>\rho x\}}{\rP\{W_1>2x\}}$. In view of (\ref{qq}),
$A(y)<\infty$ and $\lix
A(x)=\dfrac{k_\beta^n}{1-k_\beta}\left(\dfrac{2}{\rho}\right)^\beta$.
Now we have for $x\geq y$
$$\dfrac{\rP\{(\sum_{m=n+1}^\infty d_m^2)^{1/2}>2x\}}{\rP\{W_1>2x\}}\leq A(y)+\nu
\dfrac{\rP\{(\sum_{m=n+1}^\infty
d_m^2)^{1/2}>x\}}{\rP\{W_1>x\}}.$$ Iterating the latter inequality
gives that for $k=0,1,\ldots$ $$\underset{x\in [2^k y,
2^{k+1}y]}{\sup}\dfrac{\rP\{(\sum_{m=n+1}^\infty
d_m^2)^{1/2}>x\}}{\rP\{W_1>x\}}\leq $$ $$\leq
A(y)(1+\nu+\ldots+\nu^{k-1})+\nu^k \underset{x\in
[y,2y]}{\sup}\dfrac{\rP\{(\sum_{m=n+1}^\infty
d_m^2)^{1/2}>x\}}{\rP\{W_1>x\}}.$$ Let $k\to \infty$ to obtain
$$\underset{x\to\infty}{\lim \sup}\dfrac{\rP\{(\sum_{m=n+1}^\infty
d_m^2)^{1/2}>x\}}{\rP\{W_1>x\}}\leq A(y)(1-\nu)^{-1}.$$ Now
sending $y\to\infty$ gives
\begin{equation}
\label{ta} \underset{x\to\infty}{\lim
\sup}\dfrac{\rP\{(\sum_{m=n+1}^\infty
d_m^2)^{1/2}>x\}}{\rP\{W_1>x\}}\leq
\dfrac{k_\beta^n}{(1-k_\beta)(1-\nu)}\left(\dfrac{2}{\rho}\right)^\beta
=const \  k_\beta^n.
\end{equation}
For any $\lambda\in (0,1)$ and any $n\in \mn$ $$\rP\{S>x\}\leq
\rP\{S_n>(1-\lambda)^{1/2}x\}+\rP\{(\sum_{k=n+1}^\infty
d^2_k)^{1/2}>\lambda^{1/2}x\}.$$ Therefore
$$\underset{x\to\infty}{\lim \sup}\dfrac{\rP\{S>x\}}{\rP\{W_1>x\}}\overset{(\ref{ss}),
(\ref{ta})}{\leq}(1-\lambda)^{-\beta/2}\sum_{m=0}^{n-1}k_\beta^m+
const \ \lambda^{-\beta/2}k_\beta^n.$$ Let $n\to\infty$ and then
$\lambda\to 0$ to get the desired bound for the upper limit
$$\underset{x\to\infty}{\lim
\sup}\dfrac{\rP\{S>x\}}{\rP\{W_1>x\}}\leq \dfrac{1}{1-k_\beta}.$$
This completes the proof for $S$.\newline (II) From the already
proved relation
\begin{equation}
\label{ma} \rP\{W^\ast>x\}\sim (1-k_\beta)^{-1}\rP\{W_1>x\}\sim
(1-k_\beta)^{-1}x^{-\beta}L(x),
\end{equation}
it follows that $\rE W^\ast<\infty$ which in turn ensures the
uniform integrability of $W_n$.

Let us now prove (\ref{ww}). Since $W^\ast\geq W$ a.s.,
$\rE(W^\ast-x)^+\geq \rE(W-x)^+$ for any $x\geq 0$. (\ref{ma})
together with Proposition 1.5.10 \cite{BGT} implies that
$$\rE(W^\ast-x)^+ \sim
(\beta-1)^{-1}(1-k_\beta)^{-1}x\rP\{W_1>x\}.$$ Therefore,
\begin{equation}
\label{low} \underset{x\to\infty}{\lim
\sup}\dfrac{\rE(W-x)^+}{x\rP\{W_1>x\}}\leq
\dfrac{1}{(\beta-1)(1-k_\beta)}.
\end{equation}
For each $x\geq 1$ define the stopping time $\nu_x$ by
\begin{equation*}
\nu_x := \left\{%
\begin{array}{ll}
    \inf\{n \geq 1: W_n>x\}, & \hbox{if $W^\ast>x$;} \\
    +\infty, & \hbox{otherwise.} \\
\end{array}%
\right.
\end{equation*}
The random variable $W$ closes the martingale $W_n$. Hence, for
each $x\geq 1$
\begin{equation*}
\rE (W-x)1_{\{\nu_x<\infty\}}=\rE
(W_{\nu_x}-x)1_{\{\nu_x<\infty\}},
\end{equation*}
and hence $$\rE (W-x)^+ \geq
\rE(W_{\nu_x}-x)^+1_{\{\nu_x<\infty\}}.$$ We now transform the
right-hand side into a more tractable form
$$\rE(W_{\nu_x}-x)^+1_{\{\nu_x<\infty\}}=\sum_{k=1}^\infty \rE (W_k-x)^+1_{\{\nu_x=k\}}=\rE\sum_{k=1}^\infty \rE_{k-1}((W_k-x)^+
1_{\{\nu_x\geq k\}})=$$$$=\rE\sum_{k=1}^\infty 1_{\{\nu_x\geq
k\}}\rE_{k-1}(W_k-x)^+=\rE\sum_{k=1}^{\nu_x} \rE_{k-1}(W_k-x)^+.$$
From (\ref{repr}) and (\ref{basic}) with $Q$ replaced by $W_1$ it
follows that for $k=2,3,\ldots$ $$\rP_{k-1}\{W_k>y\}\sim
\sum_{|u|=k-1}Y_u^\beta \rP\{W_1>y\} \ \ \text{a.s.}$$ An appeal
to Proposition 1.5.10 \cite{BGT} gives that for $k=2,3,\ldots$
$$\rE_{k-1}(W_k-y)^+\sim (\beta-1)^{-1}\sum_{|u|=k-1}Y_u^\beta y\rP\{W_1>y\} \ \
\text{a.s.}$$ Since $\lix \nu_x=+\infty$ a.s., using Fatou's lemma
allows us to conclude that $$\underset{x\to\infty}{\lim
\inf}\dfrac{\rE(W-x)^+}{x\rP\{W_1>x\}}\geq
\underset{x\to\infty}{\lim \inf}\dfrac{\rE\sum_{k=1}^{\nu_x}
\rE_{k-1}(W_k-x)^+}{x\rP\{W_1>x\}}=$$ $$=
\dfrac{1}{\beta-1}\left(1+\sum_{k=2}^\infty \rE
\left(\sum_{|u|=k-1}Y_u^\beta\right)
\right)=\dfrac{1}{(\beta-1)(1-k_\beta)}.$$ Combining the latter
inequality with (\ref{low}) yields $$\rE(W-x)^+ \sim
(\beta-1)^{-1}(1-k_\beta)^{-1} x\rP\{W_1>x\},$$ which by the
monotone density theorem (see Theorem 1.7.2 \cite{BGT}) implies
(\ref{ww}).

The result for $M$ immediately follows from $$\rP\{W>x\}\sim
\rP\{W^\ast>x\}\sim (1-k_\beta)^{-1}\rP\{W_1>x\},$$ as $|W-1|\leq
M\leq W^\ast$ a.s. The proof of the theorem is finished.

\section{The second proof of (\ref{ww}) in the case $\beta>2$}

Assume that the assumptions of Theorem \ref{t1} hold with
$\beta>2$. By (\ref{ma}) and Theorem 1.6.5 \cite{BGT}, $$\rE
W^\ast (W^\ast-x)^+ \sim \beta
(\beta-1)^{-1}(\beta-2)^{-1}(1-k_\beta)^{-1}x^{2-\beta}L(x).$$
Since for each $x>0$ $\rE W^\ast (W^\ast-x)^+\geq \rE W(W-x)^+$
\begin{equation}
\label{le} \underset{x\to\infty}{\lim \sup} \dfrac{\rE
W(W-x)^+}{x^{2-\beta}L(x)} \leq
\dfrac{\beta}{(\beta-1)(\beta-2)(1-k_\beta)}.
\end{equation}

Lyons \cite{Lyons} constructed a probability space with
probability measure $\mq$ and proved the following equality
$$\rE_{\mq}(W_n|\mathcal{G})=1+\sum_{k=1}^n \Pi_{k-1}(S_k-1) \ \ \ \ \ \mq \ \text{a.s.},$$
where $\rE_\mq$ is expectation with respect to $\mq$; $\Pi_0:=1$,
$\Pi_k:=M_1M_2\cdots M_k, k=1,2,\ldots$; $\{(M_k, S_k):
k=1,2,\ldots\}$ are $Q$-independent copies of a random vector
$(M,S)$ whose distribution is defined by the equality
\begin{equation}
\label{eq11}\rE \sum_{|u|=1}Y_uh(Y_u, \sum_{|v|=1}Y_v)=\rE h(M,S),
\end{equation}
which is assumed to hold for any nonnegative Borel bounded
function $h(x,y)$; $\mathcal{G}$ is the $\sigma$-field that can be
explicitly described (we only note that $\sigma((M_k, S_k),
k=1,2,\ldots)\subset \mathcal{G}$). Also for any Borel function
$r$ with obvious convention when the right-hand side is infinite
or does not exist
\begin{equation}
\label{p} \rE_\mq r(W_n):=\rE W_n r(W_n) \ \ \text{and} \ \
\rE_\mq r(W):=\rE W r(W).
\end{equation}
Lyons explains his clever argument in a quite condensed form. More
details clarifying his way of reasoning can be found in
\cite{BigKypr} and \cite{IksRos}.

Since $\rP\{W_1>x\}$ regularly varies with exponent $-\beta$,
$\beta>2$ then $\rE W_1^2<\infty$. Also by (\ref{innn}) $\rE
\sum_{|u|=1}Y_u^2<1$. By Proposition 4 \cite{Iks04} the last two
inequalities together ensure that $\rE W^2<\infty$. In view of
(\ref{p}) $\rE_{\mq} W=\rE W^2$ and hence $\rE_{\mq} W<\infty$. In
Lemma 4.1 \cite{IksRos} it was proved that (\ref{inn}) holds with
$\rP$ replaced by $\mq$. This implies that $\rE_{\mq}
W^\ast<\infty$ iff $\rE_{\mq} W<\infty$. Therefore we have checked
that $\rE_{\mq} W^\ast<\infty$ which by the dominated convergence
theorem implies that
$$\rE_{\mq}(W|\mathcal{G})=1+\sum_{k=1}^\infty
\Pi_{k-1}(S_k-1)=:R \ \ \ \ \ \mq \ \text{a.s.}$$ By Jensen's
inequality, for any convex function $g$
$$\rE_{\mq}(g(W)|\mathcal{G})\geq g(\rE_{\mq}(W|\mathcal{G}))=g(R)
\ \ \ \ \ \mq \ \text{a.s.}$$ Setting $g(u):=(u-x)^+, x>0$ and
taking expectations yields
\begin{equation}
\label{tru} \rE W(W-x)^+\overset{(\ref{p})}{=}\rE_{\mq}(W-x)^+
\geq \rE_{\mq}(R-x)^+.
\end{equation}
From (\ref{eq11}) it follows that
$\rE_{\mq}M^{\beta-1}=k_\beta<1$,
$\rE_{\mq}M^{\beta-1+\epsilon}=k_{\beta+\epsilon}<\infty$ and
$$\mq\{S-1>t\}=\int_{t+1}^\infty yd\rP\{W_1\leq y\}.$$ Using the
latter equality and Theorem 1.6.5 \cite{BGT} leads to
$\mq\{S-1>t\}\sim \beta (\beta-1)^{-1}t^{1-\beta}L(t)$. Therefore
Theorem 1 \cite{Grey} applies to the perpetuity $R$ which gives
$$\mq\{R>t\}\sim
\beta
(\beta-1)^{-1}(1-\rE_{\mq}M^{\beta-1})^{-1}t^{1-\beta}L(t).$$ By
Proposition 1.5.10 \cite{BGT} $$\rE_{\mq}(R-x)^+=\int_x^\infty \mq
\{R>t\}dt \sim
\dfrac{\beta}{(\beta-1)(\beta-2)(1-k_\beta)}x^{2-\beta}L(x).$$ An
appeal to (\ref{tru}) now results in
$$ \underset{x\to\infty}{\lim \inf} \dfrac{\rE
W(W-x)^+}{x^{2-\beta}L(x)} \geq
\dfrac{\beta}{(\beta-1)(\beta-2)(1-k_\beta)}.$$ Combining this
with (\ref{le}) yields $$\rE_\mq (W-x)^+\overset{(\ref{p})}{=}\rE
W(W-x)^+\sim
\dfrac{\beta}{(\beta-1)(\beta-2)(1-k_\beta)}x^{2-\beta}L(x).$$ By
the monotone density theorem $$\mq\{W>x\}\sim
\dfrac{\beta}{(\beta-1)(1-k_\beta)}x^{1-\beta}L(x).$$ Since
$\mq\{W>x\}=\int_x^\infty yd\rP\{W\leq y\}$, integrating by parts
gives
$$\dfrac{x\rP\{W>x\}}{\mq\{W>x\}}=1-\dfrac{x}{\mq\{W>x\}}\int_x^\infty y^{-2}\mq\{W>y\}dy.$$ By Proposition 1.5.10
\cite{BGT} the right-hand side tends to $(\beta-1)\beta^{-1}$ when
$x\to\infty$. Therefore, $\rP\{W>x\}\sim (\beta-1)(\beta
x)^{-1}\mq \{W>x\}\sim (1-k_\beta)^{-1}x^{-\beta}L(x)$ as desired.
\begin{rem}
This argument seems not to work as just described when $\beta \in
(1,2]$. If $\beta\in (1,2)$ we can get a bound for the upper limit
$$ \underset{x\to\infty}{\lim \sup} \dfrac{\rE
W(W\wedge x)}{x^{2-\beta}L(x)} \leq
\dfrac{\beta}{(\beta-1)(2-\beta)(1-k_\beta)}.$$ However, we do not
know how the corresponding lower limit could be obtained. In fact,
we have not been able to find any random variable $\xi$ with
appropriate tail behaviour and such that $W\geq \xi$ in some
strong or weak sense.
\end{rem}

\section{Miscellaneous comments}

We begin this section with discussing the following problem: which
classes of point processes satisfy both (\ref{mom}) and
(\ref{var}) and which do not.

Let $h:[0,\infty)\to [0,\infty)$ be a nondecreasing and
right-continuous function with $h(+0)>0$.
\begin{ex}
Let $\{\tau_0:=0, \tau_i, i\geq 1\}$ be the renewal times of an
ordinary renewal process. In addition to the conditions on $h$
imposed above, assume that $h(0)$ is finite. Consider the point
process $\mm$ with points $\{A_i=\gamma^{-1}\log h(\tau_i),
i=1,2,\ldots\}$, where $\gamma>0$ is chosen so that
$\rE\sum_{i=1}^\infty h(\tau_i)=1$. According to Theorem 1
\cite{DoBr} $W_1=\sum_{i=1}^\infty h(\tau_i)$ has exponentially
decreasing tail. Thus while we can find $h$ and $\{\tau_i\}$ such
that (\ref{mom}) holds, (\ref{var}) always fails.
\end{ex}

The situation when $\rP\{W_1>x\}\sim x^{-\beta}L(x)$ and $\rE
W_1^\beta<\infty$ is not terribly interesting. However, if this is
the case Theorem \ref{t1} implies the one-way implication of a
well-known moment result (see \cite{Iks04} and \cite{IksNegad})
$$\rE\sum_{|u|=1}Y_u^\beta<1, \ \ \rE W_1^\beta<\infty \Leftrightarrow \rE W^\beta<\infty, \ \ \rE(W^\ast)^\beta<\infty.$$
In the subsequent examples in addition to (\ref{mom}) and
(\ref{var}) we require that $\rE W_1^\beta=\infty$. Examples
\ref{e2} and \ref{e3} essentially shows that when the number of
points in a point process is infinite, and the points are
independent or constitute an (inhomogeneous) Poisson flow, $\rE
\sum_{|u|=1}Y_u^\beta<\infty$ implies $\rE W_1^\beta<\infty,
\beta>1$.
\begin{ex}
\label{e2} Assume that $\mm$ is a point process with independent
points $\{C_i\}$, and $\rE \sum_{i=1}^\infty Y_i=1$ and for some
$\beta>1$ $\rE\sum_{i=1}^\infty Y_i^\beta<\infty$, where
$Y_i=e^{\gamma C_i}$ and $\gamma>0$. Then $\rE W_1^\beta<\infty$.

In this case $W_1=\sum_{i=1}^\infty Y_i$. Hence we must check that
$\rE (\sum_{i=1}^\infty Y_i)^\beta<\infty$. By using the
$c_\beta$-inequality let us write the (formal) inequality
\begin{equation}
\label{in} \rE(\sum_{i=1}^\infty Y_i)^\beta=\rE \sum_{i=1}^\infty
Y_i (Y_i+\sum_{k\neq i}Y_k)^{\beta-1}\leq (2^{\beta-2}\vee 1)
(\rE\sum_{i=1}^\infty Y_i^\beta+ (\rE \sum_{i=1}^\infty Y_i) \rE
(\sum_{i=1}^\infty Y_i)^{\beta-1}).
\end{equation}
For any $\beta>1$ there exists $n\in \mn$ such that $\beta\in (n,
n+1]$. We will use induction on $n$. If $\beta \in (1,2]$ then
$\rE \sum_{i=1}^\infty Y_i<\infty$ implies $\rE (\sum_{i=1}^\infty
Y_i)^{\beta-1}<\infty$. Hence by (\ref{in}) $\rE
(\sum_{i=1}^\infty Y_i)^\beta<\infty$. Assume that the conclusion
is true for $\beta \in (n, n+1]$ and prove it for $\beta \in (n+1,
n+2]$. Since $\rE\sum_{i=1}^\infty Y_i<\infty$ and
$\rE\sum_{i=1}^\infty Y_i^\beta<\infty$ we have
$\rE\sum_{i=1}^\infty Y_i^{\beta-1}<\infty$ which, by the
assumption of induction, implies $\rE(\sum_{i=1}^\infty
Y_i)^{\beta-1}<\infty$. It remains to apply (\ref{in}).
\end{ex}
\begin{ex}
\label{e3} Let $\{\tau_i, i\geq 1\}$ be the arrival times of a
Poisson process with intensity $\lambda>0$. Consider a Poisson
point process $\mm$ with points $\{B_i\}$ and assume that for any
$a\in\mr$ $\mm(a,\infty)\geq 1$ a.s. Then there exists a function
$h$ as described at the beginning of this section that
additionally satisfies $h(0)=\infty$, and $\gamma>0$ such that
$h(\tau_i)=e^{\gamma B_i}$ and $\rE W_1=\rE\sum_{i=1}^\infty
h(\tau_i)=\lambda \int_0^\infty h(u)du=1$. The distribution of
$W_1$ is infinitely divisible with zero shift and L\'{e}vy measure
$\nu$ given as follows: $\nu(dx)=\lambda h^\leftarrow (dx)$, where
$h^\leftarrow$ is a generalized inverse function. Since $\lambda
\rE\sum_{i=1}^\infty h^\beta (\tau_i)=\int_0^\infty x^\beta
\nu(dx)$, and as is well-known from the general theory of
infinitely divisible distributions, $\int_1^\infty x^\beta
\nu(dx)<\infty$ implies $\rE W_1^\beta<\infty$, we conclude that
$\rE\sum_{i=1}^\infty e^{\gamma \beta B_i}<\infty$ implies $\rE
W_1^\beta<\infty$.
\end{ex}
In the next example we point out a class of point processes which
satisfy (\ref{mom}) and (\ref{var}), and $\rE W_1^\beta=\infty$.
\begin{ex}
Let $K$ be a nonnegative integer-valued random variable with
$\rP\{K>x\}\sim x^{-\beta}L(x), \beta>1$, and $\{D_i, i\geq 1\}$
be independent identically distributed random variables which are
independent of $K$. If $\mm$ is a point process with points
$\{D_i, i=\overline{1,K}\}$ and there exists a $\gamma>0$ such
that $\rE\sum_{i=1}^K e^{\gamma D_i}=1$ and $\rE\sum_{i=1}^K
e^{\gamma \beta D_i}<1$ then according to Proposition 4.3
\cite{Sam} we have $\rP\{W_1>x\}\sim \rE e^{\gamma \beta
D_1}x^{-\beta}L(x)$.
\end{ex}
We conclude with two remarks that fit the context of the present
paper.
\begin{rem}
The tail behaviour of $\rP\{W>x\}$ and $\rP\{W^\ast>x\}$ has been
investigated in \cite{Iks04} and \cite{IksNegad}. In particular,
from those works it follows that when the condition
$\rE\sum_{|u|=1}Y_u^\beta<1$ fails, the condition
$\rP\{W_1>x\}\sim x^{-\beta}L(x)$ is not a necessary one for
either $\rP\{W>x\}\sim x^{-\beta}L(x)$ or $\rP\{W^\ast>x\}\sim
x^{-\beta}L(x)$ to hold.
\end{rem}
\begin{rem}
Let $\{X_i\}$ be the points of a point process, and let $V$ be a
random variable satisfying the following distributional equality
$$V\overset{d}{=} \sum_{i=1}^\infty X_iV_i,$$ where $V_1, V_2,
\ldots$ are conditionally on $\{X_i\}$ independent copies of $V$.
The distribution of $V$ is called a fixed point of the smoothing
transform (see \cite{Iks04} for more details, and \cite{IksJur}
and \cite{IksKim} for an interesting particular case). It is known
and can be easily checked that the distribution of $W$ is a fixed
point of the smoothing transform with
$\{X_i:i=1,2,\ldots\}=\{Y_u:|u|=1\}$. Thus (\ref{ww}) could be
reformulated as a result on the tail behaviour of the fixed points
\emph{with finite mean}. The tail behaviour of fixed points with
infinite mean deserves a special mention. Typically, their tails
regularly vary with index $\alpha\in (0,1)$, or $\int_0^x
\rP\{V>y\}dy$ slowly varies. This follows from Proposition 1(b)
\cite{Iks04} and Proposition 8.1.7\cite{BGT}.
\end{rem}

\end{document}